
\input amstex

\documentstyle{amsppt}
\redefine\aa{\alpha}
\define\bb{\beta}
\redefine\gg{\gamma}
\define\kk{\kappa}
\redefine\oA{\omega}
\redefine\oB{{\omega _1}}
\redefine\oC{{\omega _2}}
\redefine\aA{{\aleph _0}}
\redefine\aB{{\aleph _1}}
\redefine\aC{{\aleph _2}}

\define\CaB{\Bbb C(\aleph _1)}

\document

\Monograph
\NoRunningHeads
\head{Introduction}\endhead

Complete Boolean algebras proved to be an important tool in topology and set
theory. Two of the most prominent examples are $\Bbb B (\kappa ),$ the algebra
of Borel sets modulo measure zero ideal in th generalized Cantor space $\{
0,1\} ^\kappa$ equipped with product measure, and $\Bbb C(\kappa ),$ the
algebra of regular open sets in the space $\{ 0,1\} ^\kappa ,$ for $\kappa$ an
infinite cardinal. $\Bbb C(\kappa )$ is much easier to analyse than $\Bbb
B(\kappa ):$ $\Bbb C(\kappa )$ has a dense subset of size $\kappa ,$ while the
density of $\Bbb B(\kappa )$ depends on the cardinal characteristics of the
real line; and the definition of $\Bbb C(\kappa )$ is simpler. Indeed, $\Bbb
C(\kappa )$ seems to have the simplest definition among all algebras of its
size. In the Main Theorem of this paper we show that in a certain precise
sense, $\Bbb C(\aleph _1)$ has the simplest {\it structure} among all algebras
of its size, too.

\proclaim
{ Main Theorem} If ZFC is consistent then so is ZFC+$2^{\aleph _0}=\aleph
_2$+``for every complete Boolean algebra $\frak B$ of uniform density $\aleph
_1,$ $\Bbb C(\aleph _1)$ is isomorphic to a complete subalgebra of $\frak B".$
\endproclaim

There is another interpretation of the result. Let $\langle BA(\kappa
),\lessdot \rangle$ denote the class of complete Boolean algebras of uniform
density $\kappa$ quasi-ordered by complete embeddability. Then $BA(\aleph _0)$
has just one element up to isomorphism; it is $\Bbb C(\aleph _0).$ The class
$BA(\aleph _1)$ can already be immensely rich, permitting of no simple
classification; this is the case say under the continuum hypothesis. The Main
Theorem shows that the class $BA(\aleph _1)$ can have a smallest element. Note
that this smallest element must then be $\Bbb C(\aleph _1),$ since by \cite
{5, Proposition 7} $\Bbb C(\aleph _1)$ is minimal in $BA(\aleph _1).$

The techniques introduced in this paper provide us with much more information.
Most notably we get

\proclaim
{ Corollary 14} Under MA$_{\aleph _1}$ $\Bbb C(\aleph _1)$ embeds into every
complete c.c.c. Boolean algebra of uniform density $\aleph _1.$
\endproclaim

\proclaim
{ Corollary 37} Under PFA $\Bbb C(\aleph _1)$ embeds into every complete
Boolean algebra of uniform density $\aleph _1.$
\endproclaim

The search for complex objects which have to be embedded into complete Boolean
algebras of small size has been going on for some time. It has been proved
that every algebra in the class $BA(\aleph _1)$ may have to add a real \cite
{ST}, indeed a Cohen real \cite {Z}. Every uncountable Boolean algebra may
have to have an uncountable independent subset \cite {T}.

The proof of the Main Theorem is an iteration argument. The heart of the
matter lies in introducing a regular embedding of $\Bbb C(\aleph _1)$ to a
given algebra $\frak B$ of uniform density $\aleph _1$ by a sufficiently mild
forcing. This problem is solved in the first three sections. Section 1
introduces the crucial auxiliary notion of an avoidable subset of the algebra
$\frak B,$ Section 2 deals with productively c.c.c. $\frak B$ as an easier
special case, proving Corollary 14 and setting the stage for the attack at the
general case in Section 3. At the end of Section 3 we are able to demonstrate
the Main Theorem. Section 4 is devoted to a couple of relevant ZFC examples of
algebras of bigger density. Finally, Section 5 suggests several open problems.

The arguments in the paper are given a nested structure, in the order of
priority Theorem, Lemma, Claim. It is advisable for example on the first
reading of the proof of Theorem X to leave out the arguments for the Lemmas.
Our notation follows the set-theoretic standard as set forth in \cite {4}.
Throughout the paper we work with separative partially ordered sets
representing dense subsets of Boolean algebras in question rather than with
the algebras themselves. ``Algebra" stands for ``complete Boolean algebra" and
``embedding," ``embeds" stand for ``complete embedding," ``completely embeds".
In a forcing notion we write $p\geq q$ to mean that $q$ is more informative
than $p$ (i.e., the Western way); $p\perp q$ to mean that $p$ and $q$ are
incompatible, that is, no $r$ is less than both $p$ and $q.$ All partial
orders in this paper will have a maximal element by default, denoted by 1. A
poset $P$ is {\it separative} if for $p\not \leq q$ there is $r\leq p,r\perp
q.$ We say t!  hat $P$ has {\it uniform density} $\kappa$ if $|P|=\kappa$ and
for no $p\in P, R\in [P]^{<\kappa }$ $R$ is dense below $p.$ An algebra has
uniform density $\kappa$ if it has a dense subset of uniform density $\kappa
.$ If $p\in P$ then $P\restriction p$ stands for $\{ r\in P: r\leq p\} .$ We
write $P\lessdot Q$ ($P$ {\it embeds into} $Q)$ if there is $\dot H,$ a
$Q$-name such that $Q\Vdash$``$\dot H\subset \check P$ is generic over $V".$
Thus $P\lessdot Q$ iff $RO(P)$ embeds into $RO(Q)$ and we can reasonably use
$\lessdot$ for embedding of algebras. $\Bbb C(\kappa )$ is construed as
$RO(C_\kappa ),$ where $C_{\kappa}=\{ h:h$ is a function and $dom(h)\in
[\kappa ]^{<\aleph _0},\ rng(h)\subset 2\}$ ordered by reverse inclusion. For
an ordinal $\alpha$ and a set $X$ of ordinals we write $\alpha ^{*X}$ for
$min(X\setminus \alpha +1)$. $H_\kappa$ is the collection af all sets
hereditarily of size $<\kappa .$ For two models $M,N$ $M\prec N$ means that
$M$ is an elementary submod!  el of $N$ and the special predicates will be
often understood from the context.  \qed C105 marks end of proof of Claim 105,
\qed T61 marks end of proof of Theorem 61 etc. 

The results in this paper were obtained during the meeting of the two authors
at Rutgers University in September 1994 and the week following it. The second
author would like to thank Rutgers University for its hospitality during this
time. Theorem 8, Definition 20 and Lemma 21 are due to the first author, Lemma
42 is due to both authors independently and the other results are due to the
second author. The results of this paper appeared in the Chapter 2 of second
author's Ph.D. thesis.

{1. The overall strategy}

Of course, the proof of the Main Theorem is by a forcing iteration argument.
The basic challenge is, given a poset $P$ of uniform density $\aB,$ to find a
sufficiently mild forcing $Q$ such that $Q\Vdash$``$C_\aB\lessdot P".$ Then we
can hope to iterate the procedure to obtain a model for the desired statement.

The following notion plays a very important role in our argument.

\proclaim
{Definition 1} Let $P$ be an arbitrary poset. A set $D\subset P$ is called
almost avoidable if for every $p\in P$ there is a finite set $tr(p)\subset D,$
called a trace of $p$ in $D,$ such that for any $b\in [D]^{<\aA}$ with $b\cap
tr(p)=0$ there is $p^\prime \leq p$ which is incompatible with every element
of $b.$
\endproclaim

For example, any finite set $D\subset P$ is almost avoidable (set $tr(p)=D$
for every $p\in P)$ and any antichain $D\subset P$ is almost avoidable (set
$tr(p)=\{ r\},$ where $r\in D$ is some element of $D$ compatible with $p,$ for
every $p\in P).$ However, we shall be interested in finding a {\it dense}
almost avoidable set $D\subset P.$ Here is a canonical example of such a
situation. Let $P$ be the Cohen poset ${^{<\oA}\oA}$ ordered by reverse
extension. Then $P,$ as a subset of itself, is almost avoidable; just set
$tr(s)=\{ t\in P:t\subset s\} .$ If $b$ is a finite set in ${^{<\oA}\oA}$ with
$b\cap tr(s)=0$ then there is a one-step extension of the sequence $s$
avoiding every element of $b.$

The relevance of Definition 1 to our problem is explained in the following two
lemmas. They show that the statement ``a poset $P$ has a dense almost
avoidable subset" is a good approximation of ``$C_\aB\lessdot P".$

\proclaim
{Lemma 2} Let $P$ be a poset of size $\kappa$ such that $C_\kk\lessdot P.$
Then $P$ has a dense almost avoidable subset.
\endproclaim

\demo
{Proof} Let $P$ be an arbitrary poset of size $\kk,P=\{p_\aa:\aa\in\kk\}$ and
suppose that $C_\kk\lessdot P.$ Choose a $P$-name $\dot c$ such that
$P\Vdash$``$\dot c:\kk\to 2$ is $C_\kk$-generic" and fix the induced embedding
$e$ of $\Bbb C(\kk)$ to $RO(P).$ We define the set $D\subset P$ as follows:
for each $\aa\in\kk,$ we choose a condition $p_\alpha ^\prime \leq p_\aa$ and
a bit $i(\aa)\in 2$ such that $p_\aa^\prime\Vdash$``$\dot c(\aa)=i(\aa)";$ we
set $D=\{ p_\aa^\prime :\aa\in\kk\}.$

Now obviously the set $D$ is dense in $P.$ We must show that $D$ is almost
avoidable. To this aim, fix a condition $p\in P.$ Definition 1 calls for a
trace of $p$ in the set $D.$ We choose a finite function $h\in C_\kk$ with
$h\leq proj_{\Bbb C(\kk)}(p)$ and set $tr(p)=\{ p_\aa^\prime:\aa\in dom(h)\}.$

To see that the set $tr(p)$ has the required properties, let $b\subset D$ be a
finite set disjoint from $tr(p).$ So necessarily there is a finite set
$d\subset\kk$ disjoint from $dom(h)$ such that $b=\{ p_\aa^\prime :\aa\in
d\}.$ Let $k\in C_\kk$ be the function with $dom(k)=dom(h)\cup d$ and
$k(\aa)=h(\aa)$ for $\aa\in dom(h)$ and $k(\aa)=1-i(\aa)$ for $\aa\in d.$
Since $k\leq h\leq proj_{\Bbb C(\kk)}(p),$ in the poset $P$ there must be a
lower bound $p^\prime$ of the conditions $p$ and $e(k).$ By the choice of the
function $k,$ necessarily $p^\prime \perp p_\aa^\prime$ for $\aa\in d,$ and so
$p^\prime \leq p$ witnesses the statement of Definition 1 for $p,tr(p)$ and
$b.$ \qed L2
\enddemo

\proclaim
{Lemma 3} Let $P$ be a poset of uniform density $\kk$ with a dense almost
avoidable subset. Then $C_\kk\Vdash$``$C_\kk\lessdot \check P".$
\endproclaim

\remark
{Remark} We do not know about any useful strengthenings of Lemma 3; cf.
Problem 51.
\endremark

\demo
{Proof} Let $P$ be a poset of uniform density $\kk$ with a dense almost
avoidable subset $D.$ First, using the uniform density of $P$ we extract a
system of $\kk$ many disjoint maximal antichains of the set $D.$

\proclaim
{Claim 4} There is a system $\langle A_\gg:\gg\in\kk\rangle$ of pairwise
disjoint maximal antichains of the set $D.$
\endproclaim

\demo
{Proof} We fix a bookkeeping device, a bijection $e:P\times \kk\to\kk.$ By
induction on $\aa\in\kk,$ we construct a sequence $\langle
p_\aa:\aa\in\kk\rangle$ of pairwise distinct conditions in $D$ as follows.
Given $\aa\in\kk,\aa=e(p,\gg)$ and the sequence $\langle p_\bb
:\bb\in\aa\rangle,$ the condition $p_\aa$ is any condition in the set $D$
which is less than $p$ and does not appear on the sequence $\langle p_\bb
:\bb\in\aa\rangle.$ It is possible to choose such a condition since the set
$D,$ unlike the set $\{ p_\bb :\bb\in\aa\},$ is dense below the condition $p.$

By the construction, for $\gg\in\kk$ the sets $D_\gg=\{ p_\aa:\aa\in
e''P\times\{ \gg\}\}\subset D$ are pairwise disjoint dense in $P.$ The Claim
follows by choosing a maximal antichain $A_\gg\subset D_\aa$ for each
$\gg\in\kk.$ \qed C4
\enddemo

Fix a system $\langle A_\gg:\gg\in\kk\rangle$ of antichains as in Claim 4. So
we have $A_\gg\subset D$ is a maximal antichain of the poset $P$ by the
density of $D.$

\proclaim
{Definition 5} A forcing $Z$ is defined by $Z=\{ z:z$ is a function with
$dom(z)\in [\bigcup _{\gamma \in \kk}A_\gamma ]^{<\aleph _0},\ rng(z)\subset
2\} ;$ order by reverse extension.
\endproclaim

\remark
{Explanation} Essentially, we force a $P$-name for a $C_\kk$-generic sequence
$\langle \dot c_\gg:\gg \in\kk\rangle$ by finite conditions. Given
$\gg\in\kk,$ the name $\dot c_\gg$ will be a function from $A_\gg$ to $2;$ for
a condition $z\in Z,$ the function $z\restriction A_\gg$ is a finite piece of
the future $\dot c_\gg.$
\endremark

Obviously, the forcing $Z$ is isomorphic to $C_\kk.$ Thus we will have proven
the Lemma once we show that $Z\Vdash$``$C_\kk\lessdot P".$ If $H\subset Z$ is
a generic filter and $\gg\in\kk$ then $\dot c_\gamma =\bigcup H\restriction
A_\gamma$ is a $P$-name for an element of $2.$ We show that $Z\Vdash
P\Vdash$``$\langle \dot c_\gamma :\gamma \in \kk\rangle$ is
$C_{\kk}$-generic". To this end, fix $z_0\in Z,$ $z_0\Vdash$``$\dot E\subset
C_{\kk}$ is open dense" and $p_0\in P.$ We find $z_1\leq z_0,p_1\leq p_0$ so
that $z_1\Vdash _{Z}p_1\Vdash _P$``$\langle \dot c_\gamma :\gamma \in
\kk\rangle$ meets $\dot E",$ proving the Lemma. Choose a trace $tr(p_0)$ of
$p_0$ in the dense set $D\subset P$ and let $d=\{ \gamma \in \kk:A_\gamma \cap
tr(p_0)\neq 0\} ;$ thus the set $d$ is finite. For the rest of the proof we
adopt the following piece of notation: for two functions $h,k$ the symbol
$h\vec \cup k$ stands for the unique function with domain $dom(h)\cup dom(k)$
which is equal to $k$!  on $dom(k)$ and equal to $h$ on $dom(h)\setminus
dom(k).$

\proclaim
{Claim 6} There are a condition $z_{1/2}\leq z_0$ in $Z$ and $h\in C_{\kk}$
such that for any function $k:d\to 2$ we have $z_{1/2}\Vdash$``$\check h\vec
\cup \check k \in \dot E".$
\endproclaim

\demo
{Proof} Let $n=|d|$ and $\langle k_j:j\in 2^n\rangle$ enumerate $^d2.$ By
induction on $j\in 2^n+1$ we construct $w_j\in Z, h_j\in C_{\kk}$ so that
\roster
\item $w_0=z_0,h_0=0$
\item $w_j$'s are decreasing in $Z,$ $h_j$'s are decreasing in $C_{\kk}$
\item for $j\in 2^n$ we have $w_j\Vdash _{Z}$``$\check h_{j+1}\vec \cup \check
k_j\in \dot E".$ 
\endroster

There is no problem in the induction. $z_{1/2}=w_{2^n},h=h_{2^n}$ witness the
statement of the Claim. \qed C6
\enddemo

Pick $z_{1/2}\leq z_0,h\in C_{\kk}$ as in the Claim. By properties of the
trace we can find $p_{1/2}\leq p$ so that for every $r\in
dom(z_{1/2})\setminus tr(p_0)$ we have $r\perp p_{1/2}.$ We strengthen
$p_{1/2}$ to $p_1$ such that for every $\gamma \in dom(h)$ there is an element
of the antichain $A_\gamma$ above $p_1;$ denote this unique element by
$p^\gamma .$ Define a condition $w\in Z$ by $dom(w)=\{ p^\gamma :\gamma \in
dom(h)\} ,$ $w(p^\gamma )=h(\gamma )$ and set $z_1=w\vec \cup z_{1/2}.$ Thus
$z_1\in Z$ and moreover $z_1\leq z_{1/2}\leq z_0.$ The following Claim
completes the proof of the Lemma:

\proclaim
{Claim 7} $z_1\Vdash _{Z}p_1\Vdash _P$``the function $\gamma \mapsto \dot
c_\gamma ,\gamma \in dom(h)$ is in $\dot E".$ 
\endproclaim

\demo
{Proof} Comparing the function $\gamma \mapsto \dot c_\gamma ,\gamma \in
dom(h)$ to $h,$ I find that $z_1\Vdash _{Z}p_1\Vdash _P$``$\dot c_\gamma \neq
h(\gamma )$ implies $p^\gamma \in dom(z_{1/2})".$ By construction of
$p_{1/2},$ $\{ \gamma \in \omega _1:p^\gamma \in dom(z_{1/2})\} \subset d.$
Therefore $z_1\Vdash _{Z}p_1\Vdash _P$``$\dot c_\gamma = h(\gamma )$ for all
$\gamma \in dom(h)\setminus d".$ By our choice of $h$ and $z_{1/2}$ we have
$z_1\Vdash _{Z}p_1\Vdash _P$``the function $\gamma \mapsto \dot c_\gamma
,\gamma \in dom(h)$ is in $\dot E",$ i.e. the statement of the Claim. \qed
C7,L3
\enddemo 

\enddemo

This brings us back to our original task. Fix a poset $P$ of uniform density
$\aB.$ We construct a two-step iteration $Q=Q_0*C_\aB=Q_0\times C_\aB.$ The
forcing $Q_0$ serves to introduce a dense almost avoidable subset to $P.$ By
Lemma 3, we then have $Q\Vdash$``$C_\aB\lessdot \check P".$ In the next
section we show that in the special case of a productively c.c.c. poset $P,$
the most optimistic variation of the above scenario works. In section 3, we
work on the general case, which is somewhat harder and technically more
requiring.

\head
{2. Productively c.c.c. posets}
\endhead

\proclaim
{Theorem 8} Let $P$ be a separative productively c.c.c. poset with uniform
density $\aleph _1.$ Then there is a c.c.c. forcing $Q$ such that
$Q\Vdash$``$C_{\aleph _1}\lessdot P".$
\endproclaim

\demo
{Proof} Fix a productively c.c.c. separative poset $P$ of uniform density
$\aB.$ As we have seen in the previous section, we have to introduce a dense
avoidable subset to $P.$ To begin with, we stratify the poset a little. We fix
a sequence $\langle r_\alpha :\alpha \in \omega _1\rangle$ so that
\roster
\item $\{ r_\alpha :\alpha \in \omega _1\} \subset P$ is dense
\item $\forall \beta \in \alpha \in \omega _1\ r_\beta \not \leq r_\alpha$
\endroster
together with a closed unbounded set $C\subset \omega _1$ with all $\alpha \in
C$ satisfying

$$\langle \{ r_\beta :\beta \in \alpha ^{*C}\} ,\{ r_\beta :\beta \in \alpha
\} ,\leq \rangle \prec \langle \{ r_\beta :\beta \in \omega _1 \} ,\{ r_\beta
:\beta \in \alpha \} ,\leq \rangle .\eqno (P1)$$

 Let us remind the reader that for an ordinal $\nu$ and a set $X$ of ordinals,
we use the notation $\nu^{*X}=min(X\setminus(\nu+1)).$ The desired forcing $Q$
will be defined as an iteration $Q_0*\dot C_\aB$ of two c.c.c. forcings.

\proclaim
{Definition 9} $Q_0$ is the set of all functions $q$ satisfying the following:
\roster
\item"{(D9.1)}" $dom(q)\in [C]^{<\aleph _0},$ $\forall \alpha \in dom(q)\
q(\alpha )=\langle p_\alpha ^q, g_\alpha ^q\rangle ;$ if no confusion is
possible we drop the superscript $q$  
\item"{(D9.2)}" $\forall \alpha \in dom(q)\ p_\alpha \in \{ r_\beta :\alpha
\leq \beta <\alpha ^{*C}\} ,$ $g_\alpha \subset dom(q)\cap \alpha$ 
\item"{(D9.3)}" $\forall \alpha \in dom(q)\ \exists p_\alpha ^\prime \leq
p_\alpha \ \forall \beta \in (\alpha \cap dom(q))\setminus g_\alpha \ p_\alpha
^\prime \perp p_\beta .$  
\endroster
Order is by reverse extension. we set $\bar q=\{ p\in P: \exists \alpha \in
dom(q)\ p=p_\alpha \} .$
\endproclaim

\remark
{Explanation} So this is a rather straightforward try to force a dense almost
avoidable subset $D\subset P$ with finite conditions. For $q\in Q_0,$ the set
$\bar q$ is a finite piece of the future set $D.$ In the generic extension, we
will need to produce a trace of $p_\aa^q$ in $D.$ This is the role of
$g_\aa^q:$ we shall set $tr(p_\aa^q)=\{ p_\aa^q\} \cup \{ p_\bb^q:\bb\in
g_\aa^q\} .$ Note that it is enough to produce traces for a dense set of
conditions in $P.$
\endremark

\proclaim
{Lemma 10} $Q_0$ is c.c.c.
\endproclaim

\demo
{Proof} Assume for contradiction that $\{ q_\xi :\xi \in \omega _1\}$ is an
antichain in $Q_0;$ without loss of generality $|q_\xi |=n$ for all $\xi \in
\omega _1$ for some fixed $n\in \omega .$ Applying $\Delta$-system argument to
$\{ dom(q_\xi ):\xi \in \omega _1\}$ and using pigeonhole principle repeatedly
we can obtain $a\in [\omega _1]^{<\aleph _0},q\in Q_0,dom(q)=a$ and a set
$A\subset \omega _1$ of full cardinality so that for every $\xi <\nu$ in $A$
we have $q_\xi \cap q_\nu =q$ and $max(dom(q_\xi ))<min(dom(q_\nu )\setminus
a).$ Note that now no confusion is possible with the notation $p_\alpha
=p_\alpha ^{q_\xi }$ if $\alpha \in dom(q_\xi )\setminus a$ for some $\xi \in
A,$ since this $\xi$ is unique.
\proclaim
{Claim 11} For each $\xi\in A$ and each $\aa\in dom(q_\xi )\setminus a,$ there
is a condition $p_\aa^\prime \leq p_\aa$ with the following properties:
\roster
\item"{(C11.1)}" $p_\aa^\prime \leq p_\aa$ witnesses (D9.3) for $\aa$ and
$q_\xi$ 
\item"{(C11.2)}" for each $\delta \in dom(q_{\xi ^{*A}})\setminus a$ we have
$p_\aa ^\prime \perp p_\delta.$ 
\endroster
\endproclaim

\demo  
{Proof} Fix $\xi\in A$ and $\aa\in dom(q_\xi )\setminus a$ as required in the
Claim. First we choose a condition $p_\aa^{\prime 0}\leq p_\aa$ witnessing
(D9.3) for $q_\xi$ and $\aa.$ By the elementarity properties of $C$ (P1) we
can require that $p_\aa^{\prime 0}\in \{ r_\bb:\bb\in\aa^{*C}\}.$ Now let
$\delta _0<\delta _1<$ $\dots <\delta _i<\dots ,$ $i<n-|a|,$ be a list of all
ordinals in $dom(q_{\xi ^{*A}})\setminus a.$ By induction on $i\leq n-|a|$ we
build $p^{\prime i}\in P$ so that
\roster
\item $p^{\prime 0}\geq p^{\prime 1}\geq \dots$
\item $p^{\prime i}\in \{ r_\beta :\beta \in \delta _i^{*C}\}$
\item $p^{\prime i+1}\perp p_{\delta _{i+1}}$ for $i<n.$
\endroster
$p^{\prime 0}=p_\aa^{\prime 0}$ already satisfies all of (1),(2),(3). Given
$p^{\prime i},i<n-|a|,$ we can choose $p^{\prime i+1}\leq p^{\prime i}$ as
required since by (2) and the choice of $\langle r_\beta :\beta \in \omega
_1\rangle$ we have $p^{\prime i}\not \leq p_{\delta _{i+1}}.$ Notice that
$p_{\delta _{i+1}}\in \{ r_\beta :\delta _{i+1}\leq \beta < \delta
_{i+1}^{*C}\} .$ To make (2) hold for $i+1$ we use (P3) again and find
$p^{\prime i+1}\in \{ r_\beta :\delta _{i+1}\leq \beta <\delta _{i+1}^{*C}\}
.$

We set $p_\alpha ^{\prime}=p^{\prime n-|a|}.$ Thus $p_\alpha ^{\prime}\leq
p_\alpha ^{\prime 0}$ is still a witness of (D9.3) for $q_\xi$ and $\alpha$
and moreover $p_\alpha ^{\prime}\perp p_\delta$ for all $\delta \in dom(q_{\xi
^{*A}})\setminus a.$ \qed C11
\enddemo

Fix a sequence of $p_\aa^\prime$'s for $\aa\in dom(q_\xi)\setminus a,\xi \in
A$ as in the Claim. Let $B\subset A$ be a set of cardinality $\aB$ such that
for all $\xi \in B$ we have $\xi ^{*B}>\xi ^{*A}.$ For each ordinal $\xi\in
B,$ let $\langle \alpha _{i,\xi }:i<n-|a|\rangle$ be an increasing list of all
ordinals in $dom(q_\alpha )\setminus a.$ The collection $\{ \langle p_{\alpha
_{i,\xi }}^{\prime}:i<n-|a|\rangle :\xi \in B\}$ is not an antichain in
$P^{n-|a|}$ since the poset $P$ is productively c.c.c. and the collection in
question is indexed by the uncountable set $B.$ Thus we may pick ordinals $\xi
<\nu$ in $B$ so that $p_{\alpha _{i,\xi }}^{\prime}$ is compatible with
$p_{\alpha _{i,\nu }}^{\prime}$ for all $i<n-|a|.$

\proclaim
{Claim 12} The conditions $q_{\xi ^{*A}},q_\nu$ are compatible in $Q_0.$
\endproclaim

\demo
{Proof} Set $\mu =\xi ^{*A}$ and $q=q_\mu \cup q_\nu .$ We need to verify that
$q\in Q_0;$ then $q$ is the needed lower bound of $q_\mu ,q_\nu ,$ proving the
Claim. The only difficulty here is checking (D9.3) for $q.$ We split into two
cases: $\alpha \in dom(q_\mu )$ and $\alpha \in dom(q_\nu )\setminus a.$ In
the former case, $p_\alpha ^\prime$ witnessing (D9.3) for $g_\mu$ and $\alpha$
will do, since the only new values for $q$ as compared to $q_\mu$ are above
$\alpha .$ In the latter case, we find $i<n-|a|$ with $\alpha =\alpha _{i,\nu
}$ and set $p_\alpha ^{\prime \prime}$ to be a common lower bound of
$p_{\alpha _{ i,\xi }} ^{\prime}$ and $p_{\alpha _{ i,\nu }} ^{\prime},$ which
exists by the choice of $\xi <\nu .$ We claim that $p_\alpha ^{\prime \prime}
\leq p_\alpha$ witnesses (D9.3) for $q$ and $\aa:$
\roster
\item Let $\beta \in (dom (q_\nu )\cap \alpha )\setminus g_\alpha ^{q_\nu }.$
Then $p_\beta \perp p _\alpha ^\prime$ and as $p_\alpha ^{\prime \prime}\leq
p_\alpha ^\prime$ we have $p_\beta \perp p_\alpha ^{\prime \prime}$ as well.  
\item Let $\beta \in dom(q_\mu )\setminus a.$ Then by construction of
$p_{\alpha _{ i,\xi }} ^{\prime}$ (Claim 2.8) we have $p_\beta \perp p_{\alpha
_{ i,\xi }} ^{\prime}$ and as $p_\alpha ^{\prime \prime}\leq p_{\alpha _{
i,\xi }} ^{\prime}$ we conclude that $p_\beta \perp p_\alpha ^{\prime \prime}$
again. 
\endroster
All relevant $\beta$'s from the second universal quantifier in (D9.3) for $q$
and $\alpha$ have been checked. The Claim follows. \qed C12
\enddemo

By the choice of $B$ we have that $\xi ^{*A}<\nu$ and so Claim 12 stands in
direct contradiction with our assumption on $\{ q_\xi :\xi \in \omega _1\}$
being an antichain. \qed L10
\enddemo

The forcing $Q_0$ as above is actually even productively c.c.c.  since its
definition from $\langle r_\alpha :\alpha \in \omega _1\rangle$ and $C$ is
absolute, and ``productive c.c.c." of the poset $P$ is preserved under c.c.c.
forcings.

Fix a generic filter $G\subset P$ and work in $V[G].$ We define a set
$D\subset P$ by $D=\{ \bar q:q\in G\}.$

\proclaim
{Lemma 13} The set $D\subset P$ is dense almost avoidable in $P.$
\endproclaim

\demo
{Proof} As for the density of $D,$ work in $V$ for a moment. Let $q_0\in Q_0$
and $p\in P.$ Choose $\delta \in C,\ \delta >max(dom(q_0))$ so that there is
$\alpha \in \delta$ with $r_\alpha \leq p.$ By elementarity properties of $C$
(P1) there is $\beta ,$ $\delta \leq \beta <\delta ^{*C}$ with $r_\beta \leq
r_\alpha .$ We set $q_1=q_0\cup \{ \langle \delta ,\langle r_\beta
,dom(q_0)\rangle \rangle \} .$ We have that $q_1\in Q_0,q_1\leq q_0$ and
$q_1\Vdash$`` there is an element of $\dot D$ below $\check p".$ The density
of the set $D\subset P$ follows by a genericity argument.

As for the almost avoidability, let $p\in P.$ We shall produce a trace of $p$
in the set $D$ with the required properties. There is $q_0\in G$ and $\alpha
\in dom(q_0)$ such that $p_\alpha ^{q_0}\leq p.$ we claim that the trace
$tr(p)=\{ p_\aa^{q_0}\} \cup \{ p_\xi^{q_0}:\xi \in g_\alpha ^{q_0}\}$ does
the trick. To see this, choose $b\in [D]^{<\aleph _0}$ disjoint from $tr(p).$
One can find $q_1\leq q_0, q_1\in G$ with $b\subset \overline {q_1}.$ Notice
that $p_\aa^{q_0}=p_\alpha ^{q_1}$ and $g_\alpha ^{q_0}=g_\alpha ^{q_1}.$
Choose $p^\prime \leq p_\aa$ witnessing (D9.3) for $q_1, \alpha .$ By
elementarity properties of the set $C$ (P1) there is such $p^\prime$ in $\{
r_\beta :\beta \in \alpha ^{*C}\} .$ Now we repeat the process from Claim 11
to get $p^{\prime \prime}\leq p^\prime$ which is incompatible with all
$p_\delta ^{q_1},$ for $\delta\in dom(q_1)\setminus (\aa +1);$ such $p^{\prime
\prime}$ will be incompatible with all elements of $\bar q_1$ except those in
$tr(p!  ).$ It follows that $p\geq p_\aa^{q_0}\geq p^{\prime \prime}\perp r$
for all $r$ in $b.$ Therefore $p^{\prime \prime}$ witnesses the desired
property of $tr(p)$ with respect to $b.$ \qed L13
\enddemo

 Note that in $V[G],$ the poset $P$ still has uniform density $\aB.$ The
reason is that this is expressible by the first-order statement ``for no
ordinals $\aa,\bb<\oB$ the set $\{ r_\xi:\xi<\aa\}$ is dense below $r_\bb",$
whose falsity is absolute between $V$ and $V[G].$ So we can use Lemma 3,
finishing the proof of Theorem 8. The forcing we have been looking for is
$Q_0*\dot C_\aB=Q_0\times C_\aB.$ \qed T8

\enddemo

\proclaim
{Corollary 14} Under MA$_{\aleph _1}$ the algebra $\Bbb C(\aleph _1)$ embeds
into all c.c.c. algebras of uniform density $\aleph _1.$
\endproclaim

\demo
{Proof} Assume MA$_{\aleph _1}$ and choose a separative c.c.c. poset $P$ of
uniform density $\aleph _1.$ Without loss of generality the underlying set of
$P$ is $\omega _1.$ By \cite {W} the poset $P$ is $\sigma$-centered and so by
Theorem 3 there is a c.c.c. $Q$ with $Q\Vdash$``$C_{\aleph _1}\lessdot \check
P".$ Choose a large regular cardinal $\kappa$ and a model $M\prec \langle
H_\kappa ,\in ,P,Q\rangle$ with $\omega _1\subset M,$ $|M|=\aleph _1.$ The
poset $Q\cap M$ is c.c.c. and so we can use MA to get a filter $G\subset Q\cap
M$ which meets all sets in $\{ D\cap M: D\in M,D\subset Q$ dense$\} ,$ since
by elementarity all of these sets are dense in $Q\cap M.$ Let $i:M\to
\overline M$ be the transitive collapse of $M,$ $\overline G=i^{\prime \prime
}G.$ Then $i\restriction (P\cup C_{\aleph _1})=id$ and $\overline G\subset
i(Q)$ is $\overline M$-generic. By our choice of $Q$ and the elementarity of
$M$ we have $\overline M[\overline G]\models$``$i(C_{\aleph _1})=C_{\aleph
_1}\lessdot i(P)=P".$ The following Claim completes the proof of the
Corollary.

\proclaim
{Claim 15} The statement $C_{\aleph _1}\lessdot P$ is upwards absolute; that
is, if $M\subset N$ are two transitive models of rich fragments of set theory,
$\aleph _1^M=\aleph _1^N,$ $P\in M$ and $M\models$``$C_{\aleph _1}\lessdot P"$
then $N$ models the same statement.
\endproclaim

\demo
{Proof} We use an alternative characterization of regular embedding:
$C_{\aleph _1}\lessdot P$ if there is a function $e:C_{\aleph _1}\to RO(P)^+$
preserving incompatibility such that for every $p\in P$ there is $h\in
C_{\aleph _1}$ such that for any $k\in C_\aB$ with $k\leq h$ the value $e(k)$
is compatible with $p$ in $RO(P).$ So we have such $e$ in $M.$ Now $C_{\aleph
_1}^M=C_{\aleph _1}^N$ and $RO(P)^M\subset RO(P)^N$ is dense; thus properties
of $e$ survive in $N,$ showing that $N\models C_{\aleph _1}\lessdot P.$ \qed
C15, Co14
\enddemo 

\enddemo

\head
{3. The general case}
\endhead

In the case of a general poset $P,$ we cannot succeed with the scenario
outlined in the previous section. The forcing $Q$ defined there has a dense
subset of size $\aB,$ and that is just too simple to work:

\proclaim
{Lemma 16} Let $P$ be a $\sigma$-closed poset and let $Q$ be a forcing of size
$\aB$ preserving $\aB.$ Then $Q\Vdash$``$P$ is $\aA$-distributive". 
\endproclaim

\demo
{Proof} Let the posets $P,Q$ be as in the assumption of the lemma. Let
$Q\Vdash$``$\langle \dot D_i:i<\oA\rangle$ is a system of open dense subsets
of $P".$ We fix a bookkeeping device, a bijection $e:Q\times\omega\to\oB$ and
construct a descending sequence $\langle p_\aa:\aa\in \oB\rangle$ of
conditions in $P$ by induction as follows:
\roster
\item $p_0=1_P$
\item for $\alpha=\beta +1,$ where $\beta =e(q,i),$ we find a condition
$q^\prime \leq q$ in $Q$ and a condition $p\leq p_\beta$ in $P$ such that
$q^\prime\Vdash$``$\check p\in\dot D_i".$ We set $p_\alpha=p.$   
\item for $\alpha$ limit we set $p_\aa\in P$ to be any lower bound of the
chain $\langle p_\beta:\bb<\aa\rangle.$ 
\endroster

By the construction, $Q\Vdash$``$\forall i<\oA\ \exists \aa\in\oB\ \check
p_\aa\in \dot D_i".$ Since the forcing $Q$ preserves $\aB,$ we have that $Q$
forces the following:``for every $i<\oA,$ let $\aa_i\in\oB$ be the least
ordinal such that $p_{\aa_i}\in \dot D_i.$ Then $\dot \aa=sup_{i<\oA}\aa_i$ is
less than $\oB.$ Therefore $p_{\dot \aa}\in \bigcap _{i<\oA}\dot D_i$ and
$\bigcap _{i<\oA}\dot D_i\neq 0."$

The previous argument relativized to any $Q\restriction q$ and $P\restriction
p,$ where $q\in Q$ and $p\in P,$ gives the Lemma. \qed L16
\enddemo

Under the Continuum Hypothesis there exists a $\sigma$-closed poset $P$ of
size $\aB,$ and as shown in Lemma 16, the forcing $Q$ as defined in the
previous section cannot force $C_\aB\lessdot P.$ Tracing the problem, we
conclude that $Q_0,$ the first component of the forcing $Q,$ collapses $\aB.$
However, we are still able to modify the forcing $Q_0$ so that we get

\proclaim
{Theorem 17} For any separative partial order $P$ of uniform density $\aleph
_1$ there is a proper, $\omega _2$-p.i.c. forcing $Q$ such that
$Q\Vdash$``$C_{\aleph _1}\lessdot \check P".$ Moreover, if GCH holds then we
can find such $Q$ of size $\aleph _2.$
\endproclaim

Here, $\oC$-p.i.c. is one of the strong forms of $\aC$-c.c. introduced by the first author \cite {She}. It will be instrumental for iteration purposes later.

The proof strategy will be the same as for Theorem 8. Given the poset $P,$ we construct a mild forcing $Q_0$ which introduces a dense almost avoidable set $D\subset P.$ Then by Lemma 3, the forcing $Q=Q_0\times C_\aB$ will be as desired. Now our $Q_0$ will be almost the same as in the previous section, only modified by side conditions in the spirit of. Now every side conditions argument consists of three ingredients: a finite conditions construction, here supplied by the poset $Q_0$ from the previous section, coherent systems of models as in Definitions and a certain notion of transcendence as in Definition. We start with disclosing the systems of models.

Let $\kk$ be an uncountable regular cardinal and fix $\ll ,$ a well-ordering of $H_\kk .$ Also, choose one distinguished element $\Delta$ of $H_\kk.$

\proclaim
{Definition 18} We say that $\frak m$ is a coherent system of models if the following conditions are satisfied:
\roster
\item"{(D18.1)}" $\frak m$ is a function, $dom(\frak m)\in [\omega _1]^{<\aleph _0}$ and for each $\alpha \in dom(\frak m)$ the value $\frak m(\alpha )$ is a finite set of isomorphic countable submodels of $\langle H_\kappa ,\in ,\ll ,\Delta\rangle$
\item"{(D18.2)}" for each $\alpha <\beta$ both in $dom(\frak m)$ we have $\forall N\in \frak m(\alpha )\ \exists M\in \frak m(\beta )\ N\in M$
\item"{(D18.3)}" for each $\alpha <\beta$ both in $dom(\frak m)$ we have $\forall M\in \frak m(\beta )\ \exists N\in \frak m(\aa )\ N\in M.$
\endroster

We consider the set $\frak M$ of all coherent systems of models to be ordered by $\geq ,$ the reverse coordinatewise extension. That is, $\frak n\geq \frak m$ if $dom(\frak n)\subset dom(\frak m)$ and for each $\alpha \in dom(\frak n)$ I have $\frak n(\alpha )\subset \frak m(\alpha ).$
\endproclaim

The poset $\frak M$ is a subset of $H_\kk$ and it is not necessarily separative. Its definition has three parameters: the cardinal $\kk,$ the well-ordering $\ll$ and the distinguished element $\Delta.$ The following Definition is motivated by some technical considerations. For a detailed treatment, see \cite {Z}.

\proclaim
{Definition 19} Let $M\prec \langle H_\kk ,\in ,\ll,\Delta\rangle$ be a countable model and let $\frak m\in \frak M$ be such that $M\in\frak m(M\cap \omega _1).$ Then we define the following notions:
\roster
\item"{(D19.1)}" $pr_M(\frak m),$ the projection of $\frak m$ into $M\cap \frak M.$ This is the function defined by $dom(\frak n)=dom(\frak m)\cap M$ and $N\in \frak n(\alpha )$ iff there are models $N=N_0\in N_1\in \dots \in N_k=M$ such that $N_i\in \frak m(\alpha _i),$ where $\alpha =\aa _0<\aa _1<\dots <\aa _k=M\cap \oB$ is an increasing list of all ordinals in $dom(\frak m)$ between $\aa$ and $M\cap \oB.$
\item"{(D19.2)}" A system $\frak m$ is said to be $M$-full if for each $\aa\in M\cap dom(\frak m)$ and each $N\in \frak m(\aa )\setminus pr_M(\frak m)(\aa)$ there is $\overline M\in\frak m(M\cap\oB )$ such that $N\in \overline M$ and $i(N)\in pr_M(\frak m)(\aa ),$ where $i:\overline M\to M$ is the unique isomorphism of $\overline M$ and $M.$
\endroster
\endproclaim

Obviously, $pr_M(\frak m)\in \frak M\cap M.$ The idea behind this definition is that $pr_M(\frak m)$ should be a system in $M$ which grasps all the information about $\frak m$ understandable from within $M.$ 

Now here is the promised notion of transcendence over a countable submodel $M\prec H_\kappa.$

\proclaim
{Definition 20} Let $P$ be a separative partially ordered set. A set $R\subset P$ is {\it small} if for every $a\in [R]^{\aleph _0}$ there is $b\in [P]^{<\aleph _0}$ such that for every $r\in a$ there is $p\in b$ with $r\geq p.$
\endproclaim
    
Some elementary observations: principal filters in the poset $P$ are small; and a small set cannot contain an infinite antichain. A good example of a small set is a cofinal branch in a tree of height $\oB.$ Obviously, the set of all small subsets of $P$ is an ideal. The idea behind Definition is that if the poset $P$ is complicated enough, the small sets cannot capture the structure of $P.$ This is recorded in the following:

\proclaim
{Lemma 21} Assume that a poset $P$ has no countable locally dense subsets and let $\frak I$ denote the $\sigma$-ideal on $P$ generated by the small subsets of $P.$ Then for every $R\in \frak I$ the set $P\setminus R\subset P$ is dense; in other words, for every $p\in P$ the set $P\restriction p$ is $\frak I$-positive.
\endproclaim

\remark
{Remark} Say that a condition $p\in P$ is ``transcendental" over a countable model $M\prec H_\kk$ if $p\notin \bigcup \{ R\in M:R$ is a small subset of $P\}.$ Then the lemma says that there is a dense set of conditions in the poset $P$ ``transcendental" over $M,$ provided that $P$ has uniform density $\aleph_1.$
\endremark

\demo
{Proof} By contradiction. Assume that $p\in P,\ R\in \frak I,\ R=\bigcup _{i\in \omega }R_i$ and $P\restriction p\subset R,$ where the sets $R_i\subset P$ are small. To simplify the notation we assume that $p=1.$ There are two cases:
\roster
\item There is a c.c.c. forcing $Q$ such that $Q\Vdash$``$\check P$ is not c.c.c.". Choose such $Q$ and a $Q$-name $\dot A$ such that $Q\Vdash$``$\dot A\subset P$ is an uncountable antichain". As $Q$ preserves $\aB,$ we can find $q\in Q,\ i\in \omega$ so that $q\Vdash$``$\dot A\cap \check R_i$ is uncountable". So $R_i$ contains infinite antichains in a generic extension; therefore it must contain such an antichain in the ground model (the tree of finite sequences of pairwise incompatible elements of $R_i$ is ill-founded). So the set $R_i$ is not small, contradiction.
\item Otherwise. In particular, $P$ is productively c.c.c. We fix a large enough regular cardinal $\kappa$ and build a sequence $\langle \langle M_\alpha ,p_\alpha \rangle :\alpha \in \omega _1\rangle$ so that $M_\alpha \prec H_\kappa$ is a countable model, $p_\alpha \in P,\ P,R_i,\{ \langle M_\beta ,p_\beta \rangle :\beta \in \alpha \} \in M_\alpha$ and for no $r\in P\cap M_\alpha$ I have $r\leq p_\alpha .$ This is possible as $P\cap M_\alpha \subset P$ is not dense by my assumption on $P.$ Take $i\in \omega$ such that $A=\{ \alpha \in \omega _1:p_\alpha \in R_i\}$ is uncountable. Since $R_i\subset P$ is small, for each $\alpha \in \omega _1$ we can find a finite collection $\{ r_{\alpha ,j}:j<n_\alpha \} \subset P$ so that

$$ \forall \beta \in A\cap \aa\ \exists j<n_\alpha \ p_\bb \geq r_{\aa ,j}.\eqno { (P2)}$$

By elementarity we may and will assume that $\{ r_{\alpha ,j}:j<n_\alpha \} \subset M_\alpha .$ From the construction of $p_\aa$'s I can then conclude that $p_{\alpha ^{*A}}\not \geq r_{\alpha ,j}$ for $\alpha \in \omega _1,j<n_\alpha .$ By separativity we can strengthen all $r_{\alpha ,j}$ so that they are incompatible with $p_{\alpha ^{*A}}.$ This preserves the property (P2) of the system $\{ r_{\alpha ,j}:j<n_\alpha ,\alpha \in \omega _1 \}$ even though now $r_{\alpha ,j}$ may be outside $M_\alpha .$ Fix $n\in \omega$ and an uncountable set $B\subset \omega _1$ so that for all $\alpha \in B$ we have $n_\alpha =n$ and $\alpha ^{*B}>\alpha ^{*A}.$ Remembering that the poset $P$ is assumed to be productively c.c.c., the collection $\{ \langle r_{\alpha ,j}:j<n\rangle :\alpha \in B  \} \subset P^n$ is not an antichain in $P^n$ and we can choose $\xi<\nu$ in $B$ with $r_{\xi ,j},r_{\nu ,j}$ compatible for all $j<n.$ By (P2) there is $j<n$ so that $p_{\xi ^{*A}}\geq r_{\nu ,j}.!
 $ However, $r_{\xi ,j}$ is both incompatible with $p_{\xi ^{*A}}$ and compatible with $r_{\nu ,j},$ contradiction. \qed 
\endroster
\enddemo

Finally, we are ready to define the forcing $Q_0$ introducing a dense almost avoidable subset to a given poset $P.$ Fix a poset $P$ of uniform density $\aB.$ Without loss of generality we may assume that the universe of $P$ is $\omega_1.$ Furthermore, set $\kk=\omega_2,\Delta=P$ and fix a wellordering $\ll$ of $H_\kk.$ Below, the set $\frak M$ of systems of models will be computed using these three parameters.

\proclaim
{Definition 22} A forcing notion $Q_0$ is defined as the set of all pairs $\langle q,\frak m\rangle$ for which
\roster
\item"{(D22.1)}" $q$ and $\frak m$ are finite functions with the same domain, which is a finite subset of $\oB$
\item"{(D22.2)}" for every $\aa\in dom(q)$ the value $q(\aa)$ is a pair $\langle p_\aa^q,g_\aa^q\rangle$ where if no confusion can result, we drop the superscript $q$
\item"{(D22.3)}" $\forall \aa\in dom(q)\ p_\aa\in P$ and $g_\aa\subset dom(q)\cap \aa$
\item"{(D22.4)}" $\forall \aa\in dom(q)\ \exists p_\aa^\prime \leq p_\aa\ \forall \bb\in (dom(q)\cap\aa)\setminus g_\aa\ p_\aa^\prime \perp p_\bb$
\item"{(D22.5)}" $\frak m$ is a coherent system of models, i.e. $\frak m\in \frak M$
\item"{(D22.6)}" for every $\aa<\bb$ both in $dom(q)=dom(\frak m)$ and for every $N\in \frak m(\beta)$ we have $p_\aa\in N$
\item"{(D22.7)}" for every $\aa\in dom(q),$ for every $N\in \frak m(\aa)$ and for every small set $R\subset P$ in $N,$ we have $p_\aa\notin R.$
\endroster

The order is defined by $\langle q_0,\frak m_0\rangle \geq \langle q_1,\frak m_1\rangle$ if $q_0\subset q_1$ and $\frak m_0\geq _{\frak M}\frak m_1.$ For a condition $\langle q,\frak m\rangle \in Q_0,$ we set $\bar q=\{ p\in P:\exists \aa\in dom(q)\ p=p_\aa^q\}.$
\endproclaim

\remark
{Explanation} This may look complicated but in fact it is not. In a condition $\langle q,\frak m\rangle,$ the $q$ part is exactly like an element of $Q_0$ in the previous section, except that it ignores any stratification of the poset $P.$ The properties (D22.2,3,4) describe just this fact. The system $\frak m$ is just the control device described in Definition. Here it is tied to $q$ by (D22.6,7). The transcendence requirement (D22.7) is the main technical point in the construction.

As it was the case in the previous section, the forcing $Q_0$ serves to introduce a dense almost avoidable subset $D$ to $P.$ The set $\bar q$ is a finite piece of the future set $D,$ and the trace of $p_\aa$ will be obtained as $tr(p_\aa^q)=\{ p_\aa^q\} \cup \{ p_\bb^q:\bb\in g_\aa^q\} .$
\endremark

We start with a simple preliminary Lemma.

\proclaim
{Lemma 23} If $\langle q,\frak m\rangle\in Q_0$ and $M\prec\langle H_\kk,\in,\ll,P\rangle$ are such that $M\in\frak m(M\cap\oB)$ then there is a condition $\langle q,\frak n\rangle\in Q_0$ such that $\langle q,\frak n\rangle\leq\langle q,\frak m\rangle$ and $\frak n$ is $M$-full.
\endproclaim

\demo
{Proof}  Fix a condition $\langle q,\frak m\rangle\in Q_0.$ The $M$-full system $\frak n\leq \frak m$ will be built so that $dom(\frak m)=dom(\frak n).$ We shall start with $\frak m;$ then we gradually add some new models to the values $\frak m(\alpha ),\aa\in dom(\frak m)\cap M,$ preserving properties (D18.1,2,3), (D22.6,7) at each step. After finitely many steps, an $M$-full system $\frak n\leq \frak m$ will emerge.

Let $\aa_0<\aa_1<\dots<\aa_k=M\cap\oB$ be an increasing list of all ordinals in $dom(\frak m)$ below $M\cap\oB$ inclusive. Let $N\in \frak m(\aa_j)\setminus pr_M(\frak m)(\aa_j)$ be a model, for some $j<k.$ Then by using (D18.3) repeatedly, we can find an $\in$-chain $N_0\in N_1\in\dots\in N_j=N\in N_{j+1}\in\dots\in N_k$ such that $N_l\in \frak m(\aa_l),$ all $l\leq k.$ Let $i:N_k\to M$ be the isomorphism. We throw all models $i(N_l)$ into $\frak n(\aa_l),$ for $l<k.$ It is readily checked that this operation preserves properties (D18.1,2,3),(D22.6,7); for example, $i(N_l)$ is isomorphic to $N_l$ via $i\restriction N_l$ and if $l_1<l_2$ then $i(N_{l_1})\in i(N_{l_2}).$ We repeat this procedure for all models $N\in \frak m(\aa_j)\setminus pr_M(\frak m)(\aa_j).$ The reader can check that the resulting system $\frak n$ is as required. \qed L23

\enddemo

Now I am ready to go right into the eye of the storm. The following proof is much like some arguments in \cite {T}.

\proclaim
{Lemma 24} $Q_0$ is proper.
\endproclaim

\demo
{Proof} Choose a large regular cardinal $\lambda,$ a condition $\langle q_0,\frak m_0\rangle\in Q_0$ and a countable submodel $M\prec H_\lambda$  with $q_0,\frak m_0,\kappa, \ll ,P$ in $M.$ We shall produce a master condition $\langle q_1,\frak m_1\rangle\leq\langle q_0,\frak m_0\rangle$ for the model $M.$ Find $p\in P\setminus \bigcup \{ R\in M:R\subset P$ small $\} ;$ there is a dense set in $P$ of these due to Lemma 21. We define $q_1=q_0\cup \{ \langle M\cap \omega _1,\langle p,dom(q_0)\rangle \rangle \}$ and $\frak m_1=\frak m_0\cup \{ \langle M\cap \omega _1, \{ M\cap H_\kappa \}\rangle \}.$ 

\proclaim
{Claim 25} $\langle q_1,\frak m_1\rangle\in Q_0,\langle q_1,\frak m_1\rangle\leq\langle q_0,\frak m_0\rangle$ \qed C25
\endproclaim

We must verify that $\langle q_1,\frak m_1\rangle$ is a master condition for the model $M.$ So for any maximal antichain $A$ of $Q_0$ in $M,$ the set $A\cap M$ should be predense below $\langle q_1,\frak m_1\rangle.$ To prove this, let $A\in M$ be a maximal antichain of $Q_0$ and choose a condition  $\langle q_2,\frak m_2\rangle\leq\langle q_1,\frak m_1\rangle.$ By eventually strengthening the condition, we can assume that there is an element $x$ of $A$ above it and $\frak m_2$ is $M\cap H_\kk$-full (Lemma 23). I shall show that the element $x$ belongs actually to $A\cap M,$ finishing the proof of properness.We define a condition $\langle q_3,\frak m_3\rangle\geq\langle q_2,\frak m_2\rangle,$ a sort of projection of $\langle q_2,\frak m_2\rangle$ to the model $M.$ So, let $q_3=q_2\restriction M$ and $\frak m_3=pr_{M\cap H_\kk}\frak m_2.$

\proclaim
{Claim 26} $\langle q_3,\frak m_3\rangle \in M\cap Q_0,\ \langle q_2,\frak m_2\rangle\leq \langle q_3,\frak m_3\rangle.$ \qed C26
\endproclaim

The task now is to carefully extend the condition $\langle q_3,\frak m_3\rangle$ within $M$ to $\langle q_4,\frak m_4\rangle$ which has an element of $A$ above it and is still compatible with $\langle q_2,\frak m_2\rangle.$ Let $\alpha _0<\alpha _1 <\dots <\alpha _n$ be an increasing list of $dom(q_2)\setminus dom(q_3);$ thus $\alpha _0=M\cap \omega _1.$ For $0\leq i\leq n$ we put $p_{\alpha _i}=p_{\alpha _i}^{q_2}.$

\proclaim
{Definition 27} For all $x\in [P]^{<\aleph _0}$ simultaneously by induction on $i\in \omega$ we define sets $x^{(i)}\subset P:$
\roster
\item"{(D27.1)}" $x^{(0)}=\{ p\in P: \exists \langle q_4,\frak m_4\rangle\leq \langle q_3,\frak m_3\rangle\ \bar q_4 =\bar q_3\cup x\cup \{ p\}$ and there is an element of $A$ above $\langle q_4,\frak m_4\rangle\} .$
\item"{(D27.2)}" $x^{(i+1)}=\{ p\in P:(x\cup \{ p\} )^{(i)}$ is not small $\} .$
\endroster
\endproclaim

Notice that the collection $\{ x^{(i)}:x\in [P]^{<\aleph _0},i\in \omega \}$ is in $M\cap H_\kappa .$

\proclaim
{Claim 28} The set $0^{(n)}$ is not small in the poset $P.$
\endproclaim

\demo
{Proof} By contradiction. Assume the set is small. By induction on $0\leq i\leq n$ we prove that 
\roster
\item $Z_i=\{ p_{\alpha _0},\ p_{\alpha _1},\ \dots ,p_{\alpha _j},\ \dots ,\ j<i\} ^{(n-i)}$ is small in $P$
\item $p_{\alpha _i}\notin \{ p_{\alpha _0},p_{\alpha _1},\dots ,p_{\alpha _j},\dots ,j<i\} ^{(n-i)},$ 
\endroster
which will be a contradiction for $i=n,$ as $p_{\alpha _n}\in \{ p_{\alpha _0},p_{\alpha _1},\dots ,p_{\alpha _j},\dots ,j<n\} ^{(0)},$ as witnessed by $\langle q_2,\frak m_2\rangle.$ Now for $i=0$ we have $0^{(n)}$ is small by the assumption and $p_{\alpha _0}=p\notin 0^{(n)}\in M\cap H_\kappa$ by the choice of $p.$ Now given (1) and (2) for $i<n,$ by (D27.2) we immediately get that the set $Z_{i+1}$ is small, i.e. (1) for $i+1.$ Now by (D18.2) for the system $\frak m_2$ we find a model $N\in\frak m_2(\aa_{i+1})$ with $M\cap H_\kk\in N.$ Then $Z_{i+1}\in N$ and by (D22.7) $p_{\alpha _{i+1}}\notin Z_{i+1},$ i.e. (2) for $i+1.$ This completes the argument. \qed C28
\enddemo

We proceed with the construction of $\langle q_4,\frak m_4\rangle.$ For $0\leq j\leq n$ fix $p_{\alpha _j}^\prime \leq p_{\alpha _j}$ witnessing (D22.4) for $q_2.$ By induction on $0\leq i\leq n$ we build $r_i, p_{\alpha _j}^{\prime i},0\leq j\leq n$ so that
\roster
\item $r_i\in P\cap M, p_{\alpha _j}\geq p_{\alpha _j}^\prime \geq p_{\alpha _j}^{\prime 0} \geq p_{\alpha _j}^{\prime 1} \geq \dots \geq p_{\alpha _j}^{\prime i}$ is a decreasing sequence of elements of $P$ for all $0\leq j\leq n$
\item $r_i\in \{ r_0,r_1,\dots ,r_k, \dots ,k<i\} ^{(n-i)}\cap M$ and $Z_i=\{ r_0,r_1,\dots ,r_k, \dots ,k<i\} ^{(n-i)}$ is not small
\item $p_{\alpha _j}^{\prime i}\perp r_i$ for all $0\leq j\leq n.$
\endroster
To construct $r_0$ recall that the set $0^{(n)}$ is not small. Thus there is $a\in [0^{(n)}]^{\aleph _0}\cap M$ witnessing that. We have $a\subset M\cap H_\kappa$ and as $ \{ p_{\alpha _j}^\prime :0\leq j\leq n\}$ does not bound all elements of $a$ we can choose $r_0\in a$ with $r_0\not \geq p_{\alpha _j}^\prime$ for all $0\leq j\leq n.$ By the separativity of the poset $P$ there are $p_{\alpha _j}^\prime \geq p_{\alpha _j}^{\prime 0}$ with $r_0\perp p_{\alpha _j}^{\prime 0}.$ By (D27.2) the set $\{ r_0\} ^{(n-1)}$ is not small. The induction step from $i<n$ to $i+1$ is carried out similarly with $p_{\alpha _j}^{\prime i}$ in place of  $p_{\alpha _j}^{\prime }$ and $Z_i$ in place of $0^{(n)}.$

The induction having been carried out up to $n$ we have $r_n\in \{ r_0,r_1,\dots ,r_k, \dots ,$ $k<n\} ^{(0)}$ and so by (D27.1) applied in $M$ there exists a condition $\langle q_4,\frak m_4\rangle\in M$ such that $\langle q_4,\frak m_4\rangle\geq \langle q_3,\frak m_3\rangle,$ $\bar q_4=\bar q_3\cup \{ r_i:0\leq i\leq n\}$ and there is an element of $A\cap M$ above it.

\proclaim
{Claim 32} The conditions $\langle q_2,\frak m_2\rangle$ and $\langle q_4,\frak m_4\rangle$ are compatible.
\endproclaim

\demo
{Proof} We shall produce a lower bound of $\langle q_5,\frak m_5\rangle\in Q_0$ of the two conditions. First, we define $q_5=q_2\cup q_4.$ It is easy to see that $q_5$ is a function satisfying (D22.2,3) and such that $q_5\restriction M=q_4.$ We must check the property (D22.4). There are two cases:
\roster
\item $\aa\in dom(q_5)\cap M$ (i. e. $\aa\in dom(q_4)).$ In this case, the element $p_\aa^\prime \leq p_\aa$ witnessing (D22.4) for $q_4$ will do even for $q_5,$ since $q_5\restriction (\aa+1)=q_4\restriction (\aa+1).$
\item $\aa\in dom(q_2)\setminus M$ (i. e. $\aa\geq M\cap\oB$ and $\aa\in dom(q_2)).$ Then $\aa=\aa_j$ for some $j\leq n.$ We claim that $p_\aa^{\prime \prime}=p_{\aa_j}^\prime \leq p_\aa$ witnesses (D22.4) for $\aa$ and $q_5.$ To see this, choose an ordinal $\bb$ in $(dom(q_5)\cap \aa)\setminus g_\aa.$ Only two things can happen here. Either $\bb\in dom(q_4).$ In this case already $p_{\aa_j}^\prime \leq p_\aa$ as fixed above is incompatible with $p_\bb;$ since $p_\aa^{\prime \prime}\leq p_{\aa_j}^\prime$ we then have $p_\aa^{\prime\prime}\perp p_\bb$ as well. Or, $\bb\in dom(q_4)\setminus dom(q_2).$ Then by the above construction, $p_\bb=r_i$ for some $i\leq n$ and consequently $p_{\aa_j}^{\prime i}\leq p_\aa$ is incompatible with $p_\bb=r_i.$ Since $p_\aa^{\prime\prime}\leq p_{\aa_j}^{\prime i}$ we have $p_\aa^{\prime\prime}\perp p_\bb$ as well. All relevant $\bb$'s in the second universal quantifier in (D22.4) have been checked and (D22.4) follows.
\endroster
We still have to define the system $\frak m_5.$ Here is the place where we use the $M\cap H_\kk$-fullness of the system $\frak m_2.$ We shall have $dom(\frak m_5)=dom(\frak m_2)\cup dom(\frak m_4);$ the description of the values $\frak m_5(\aa)$ splits into two cases:
\roster
\item if $\aa\in dom(\frak m_5)\setminus M$ (i. e. $\aa\geq M\cap \oB$ and $\aa\in dom(\frak m_2))$ then $\frak m_5(\aa )=\frak m_2(\aa )$
\item if $\aa\in dom(\frak m_5)\cap M$ (i. e. $\aa\in dom(\frak m_4))$ then $\frak m_5(\aa )=\frak m_4(\aa)\cup\{ i(N):N\in \frak m_2(\aa)$ and $i:M\to\overline M$ is an isomorphism with $\overline M\in \frak m_2(M\cap \oB)\}.$
\endroster
The reader can verify that the system $\frak m_5$ is in $\frak M$ and satisfies the conditions (D22.6,7). The $M\cap H_\kk$-fullness of the system $\frak m_2$ together with our construction of $\frak m_5$ ensures that $\frak m_5\leq \frak m_2,\frak m_4$ as desired.
\qed C29
\enddemo

Now since $A\subset Q_0$ is an antichain and the conditions $\langle q_2,\frak m_2\rangle$ and $\langle q_4,\frak m_4\rangle$ are compatible, the elements of $A$ above $\langle q_2,\frak m_2\rangle$ and $\langle q_4,\frak m_4\rangle$ must be identical. However, the unique element of $A$ above $\langle q_4,\frak m_4\rangle\in M$ is in $M$ by elementarity, and I have finished the proof of Lemma 24. \qed L24
\enddemo

\proclaim
{Lemma 30} $Q_0$ has $\oC$-p.i.c.
\endproclaim

Let me remind the reader what this is all about.

\proclaim
{Definition 31} \cite {She,Ch.VIII,\S 2} A forcing $Q$ has $\omega _2$-p.i.c. 
(properness isomorphism condition) if for all large regular cardinals $\lambda$ and every $\Delta \in H_\lambda ,$ 
{\bf for every} $\gg<\delta<\omega _2,$ $q_0, h,$ $M_\gg,M_\delta$ countable 
submodels of $\langle H_\lambda ,\in ,\Delta \rangle$ with 
$\gg\in M_\gg,$ $\delta\in M_\delta,$ $Q\in M_\gg\cap M_\delta ,$ $M_\gg\cap \gg=M_\delta\cap \delta,$ 
$M_\gg\cap \omega _2\subset \delta,$ $q\in Q\cap M_\gg,$ $i:M_\gg\to M_\delta$ an 
isomorphism which is identity on $M_\gg\cap M_\delta$ {\bf there is} 
$q_1 \leq q_0,$ a master condition for $M_\gg$ such that 
$q_1 \Vdash$``$i^{\prime \prime }(\check M_\gg\cap \dot G)
=\check M_\delta\cap \dot G".$ A condition $q_1$ as above is called a symmetric master condition for $M_\gg,M_\delta.$
\endproclaim

Intuitively, we want the isomorphism $i$ to extend in the generic extension to an isomorphism $\hat i:M_\gg[G]\to M_\delta[G]$ in the most natural way: we want to set $\hat i(\tau/G)=(i(\tau))/G.$ The condition $q_1$ forces that this will indeed be an isomorphism. Perhaps at least one rather trivial example is in order: any proper forcing $Q$ of cardinality $\aB$ has $\oC$-p.i.c. This is because in any two models as in the Definition, we obtain $i\restriction Q\cap M_\gg=id.$ Therefore, every master condition $q_1$ for the model $M_\gg$ will have the required ``symmetricity" property. 

The point in such a strange property of the forcing $Q$ is that granted the Continuum Hypothesis, an $\oC$-p.i.c. forcing $Q$ has $\aC$-c.c. and preserves the Continuum Hypothesis. In fact, this is even true for short iterations of $\oC$-p.i.c. forcings:

\proclaim
{Fact 32} \cite {She,Ch.VIII,\S 2} Assume CH. If $\langle P_\alpha :\alpha \leq \omega _2,\dot Q_\alpha :\alpha < \omega _2\rangle$ is a countable support iteration of forcings such that for each $\alpha <\omega _2$ we have $P_\alpha \Vdash$``$\dot Q_\alpha$ has $\omega _2$-p.i.c." then
\roster
\item"{(F32.1)}" $\forall \alpha \in \omega _2\ P_\alpha$ has $\omega _2$-p.i.c. and $P_\alpha \Vdash$``CH"
\item"{(F32.2)}" $P_{\omega _2}$ has $\aleph _2$-c.c.
\endroster
\endproclaim

\demo
{Proof of Lemma 30} We show a little more general statement than that of Definition 31. Choose a large regular cardinal $\lambda,$ a condition $\langle q_0,\frak m_0\rangle\in Q_0$ and two isomorphic countable submodels $M_0,M_1\prec \langle H_\lambda ,\in,P\rangle$ such that $Q_0,\ll,\kk$ are in both of them and $\langle q_0,\frak m_0\rangle\in M_0.$ Let $i:M_0\to M_1$ be an isomorphism, $i(P)=P,$ $i(\ll)=\ll.$ We shall produce the desired symmetric master condition $\langle q_1,\frak m_1\rangle \leq\langle q_0,\frak m_0\rangle$ for the two models.

First, we pick $p\in P$ which does not belong to any small subset of $P$ which is in $M_0\cup M_1.$ There is a dense set of these due to Lemma 21. Now we set $q_1=q_0\cup\{ M_0\cap \oB,\langle p,dom(q_0)\rangle\rangle\}.$ We construct $\frak m_1$ as the unique function such that $dom(\frak m_1)=dom(\frak m_0)\cup\{ M_0\cap\oB\}$ and the values are defined as follows: for $\aa\in dom(\frak m_0)$ we set $\frak m_1(\aa)=\frak m_0(\aa)\cup i(\frak m_0(\aa))$ and for $\aa=M_0\cap \oB$ we set $\frak m_1(\aa)=\{ M_0\cap H_\kk,M_1\cap H_\kk\}.$ The following is immediate:

\proclaim
{Claim 33} $\langle q_1,\frak m_1\rangle\in Q_0,\langle q_1,\frak m_1\rangle\leq\langle q_0,\frak m_0\rangle$ \qed C2.30
\endproclaim

We claim that $\langle q_1,\frak m_1\rangle$ is the desired symmetric master condition for $M_0,M_1.$ Obviously, $\langle q_1,\frak m_1\rangle$ is a master condition for $M_0$ since it is stronger than the master condition described in Lemma 24. We must verify that  $\langle q_1,\frak m_1\rangle\Vdash$``$i^{\prime \prime }(\check M_0\cap \dot G)=\check M_1\cap \dot G".$ I prove that  $\langle q_1,\frak m_1\rangle\Vdash$``$i^{\prime \prime} (\check M_0\cap \dot G)\subset\check M_1\cap \dot G";$ the proof of the opposite inclusion is parallel.  
So let $x\in M_0\cap Q_0$ and let $\langle q_2,\frak m_2\rangle\leq\langle q_1,\frak m_1\rangle$ be a condition such that $\langle q_2,\frak m_2\rangle\Vdash$``$\check x\in \dot G".$ We shall obtain a condition $\langle q_2,\frak m_3\rangle$ such that $\langle q_3,\frak m_3\rangle\Vdash$``$i(\check x)\in \dot G".$ By a genericity argument, this will complete the proof. Now by eventually strengthening the condition $\langle q_2,\frak m_2\rangle$ we can assume that $\langle q_2,\frak m_2\rangle\leq x.$

\proclaim
{Claim 34} $M_0\models$``there is a system $\frak n$ such that $\frak n\leq pr_{M\cap H_\kk}\frak m_2$ and $\langle q_2\restriction M_0,\frak n\rangle\leq x".$
\endproclaim

\demo
{Proof} Notice first that the parameters of the formula--the system $pr_{M\cap H_\kk}\frak m_2,$ the condition $x$ and the finite function $q_2\restriction M_0$--are all in the model $M_0.$ The claim follows from the elementarity of $M_0,$ since the formula is witnessed in $H_\lambda$ by the system $\frak m_2\restriction M_0.$ \qed C34
\enddemo

Now let a system $\frak n$ be as in the Claim. We define a system $\frak m_3$ by $dom(\frak m_3)=dom(\frak m_2),$ for $\aa\in dom(\frak m_2)\cap M_0$ I set $\frak m_3(\aa)=\frak m_2(\aa)\cup i(\frak n(\aa))$ and for $\aa\in dom(\frak m_2)\setminus M_0$ we set $\frak m_3(\aa)=\frak m_2(\aa).$

\proclaim
{Claim 35} $\frak m_3$ is a coherent system of models and $\langle q_2,\frak m_3\rangle\in Q_0.$ \qed C35
\endproclaim

We claim that $\langle q_2,\frak m_3\rangle$ is the desired condition. First, obviously $\langle q_2,\frak m_3\rangle\leq\langle q_2,\frak m_2\rangle.$ Second, we have $\langle q_2,\frak m_3\rangle\leq i(x):$ this is because $i(q_2\restriction M_0)=q_2\restriction M_0$ and so $\langle q_2,\frak m_2\rangle\leq i(\langle q_2\restriction M_0,\frak n\rangle\leq i(x)$ by the isomorphism properties of $i.$ \qed L30 
\enddemo

Now I proceed exactly as in the previous Section.
Choose a generic filter $G\subset Q_0$ and in $V[G],$ define a set $D\subset P$ by $D=\bigcup \{ \bar q:\langle q,\frak m\rangle\in G,$ for some coherent system $\frak m\}.$

\proclaim
{Lemma 36} The set $D\subset P$ is dense almost avoidable in $P.$
\endproclaim

\demo
{Proof} This is almost exactly the same as Lemma 2.10. We show why the set $D$ is dense in the poset $P.$ Fix $p\in P$ and a condition $\langle q_0,\frak m_0\rangle\in Q_0.$ We shall produce a condition $\langle q_1,\frak m_1\rangle\leq\langle q_0,\frak m_0\rangle$ such that $\langle q_1,\frak m_1\rangle\Vdash$``there is an element of $\dot D$ below $\check p".$ The density of $D$ will then follow from a genericity argument. 
So we choose a large regular cardinal $\lambda$ and a countable elementary submodel $M\prec H_\lambda$ such that $Q,P,p,\langle q_0,\frak m_0\rangle$ are all in $M.$ By Lemma 21, there is $p^\prime \leq p$ in the poset $P$ such that $p^\prime$ does not belong to any small subset of $P$ which is in the model $M.$ We set $q_1=q_0\cup\{ \langle M\cap\oB,\{ p^\prime ,dom(q_0)\}\rangle\}$ and $\frak m_1=\frak m_0\cup\{ \langle M\cap \oB,\{ M\cap H_\kk\}\rangle\}.$ Obviously, the condition $\langle q_1,\frak m_1\rangle$ is as desired. \qed L36  
\enddemo

The proof of Theorem 17 is now finished as in the previous section with $Q=Q_0*\dot C_\aB.$ We must verify that the forcing $Q=Q_0*\dot C_\aB$ has the required properties. As in Lemma 3 $Q\Vdash$``$C_\aB\lessdot P".$ The forcing $Q$ is proper $\oC$-p.i.c. since it is an iteration of two such forcings. The last thing to check is the size of $Q.$ The forcing $Q$ as an iteration may be large, but it has a dense subset isomorphic to $Q_0\times C_\aB.$ Now under GCH, I have $|H_\aC|=\aC$ and so $|Q_0|=\aC.$ As a result, the forcing $Q$ has a dense subset of size $\aC\cdot\aB=\aC.$ Theorem 17 has been proven.\qed T17

\proclaim
{Corollary 37} Under the Proper Forcing Axiom, every complete Boolean algebra of uniform density $\aB$ contains a complete subalgebra isomorphic to $\CaB.$
\endproclaim

To simplify the proof of this, we first prove the following multipurpose Lemma.

\proclaim
{Lemma 38} The Proper Forcing Axiom implies that for every proper forcing notion $Q,$ a regular large enough cardinal $\kappa$ and a distinguished element $\Delta\in H_\kk$ there are a model $M$ so that $M\prec H_\kappa ,$ $\omega _1\subset M,\ Q,\Delta\in M$ and a filter $G\subset M\cap Q$ which is $M$-generic over $Q.$ That is, for every dense set $D\subset Q$ which is in $M,$ we obtain $G\cap D\neq 0.$
\endproclaim

\remark 
{Remark} A similar statement for MA$_{\aleph _1}$ and c.c.c. forcings is virtually trivial, since c.c.c.-ness of $Q$ is inherited by $Q\cap M:$ first, choose a model $M$ of cardinality $\aB$ and then apply MA$_{\aleph _1}$ to $Q\cap M$ and all the dense sets of $Q$ in $M.$ However, properness is not usually inherited to arbitrary subposets and we need an additional twist to complete the argument.
\endremark

\demo
{Proof} Choose a proper forcing $Q,$ a large regular cardinal $\kk$ and $\Delta\in H_\kk.$ There is a function $f:H_\kk^{<\oA}\to H_\kk$ such that if a set $M\subset H_\kk$ is closed under $f,$ then $M$ is already a submodel of $\langle H_\kk,\in,\Delta ,Q\rangle .$ So let me choose such a function $f.$ By induction on $n\in\oA$ we define simultaneously for all sets $a\in [H_\kk]^{<\aA}$ and all conditions $q\in Q$ the following finite sets $a^{(n,q)}\subset H_\kk:$
\roster
\item $a^{(0,q)}=a.$
\item The induction step from $n$ to $n+1$ is conducted as follows: I set $b=a^{(n,q)}\cup \{ x\in Q:$ there is a maximal antichain $A\subset Q$ such that $A\in a^{(n,q)},x\in A$ and $q\leq x\}.$ Then we define $a^{(n+1,q)}=b\cup f^{\prime\prime}b^{<\oA}.$
\endroster
For an integer $n$ and a set $a\in [H_\kk]^{<\aA}$ we define a subset $D_{n,a}$ of the forcing $Q$ by $D_{n,a}=\{ q\in Q:$ for every $i<n$ and every maximal antichain $A\subset Q$ with $A\in a^{(i,q)}$ there is $x\in A$ such that $q\leq x\}.$

\proclaim
{Claim 39} The sets $D_{n,a}\subset Q$ are open dense in $Q.$
\endproclaim

\demo
{Proof} The openness of $D_{n,a}$ follows straight from its definition. Note that if $q\in D_{n,a}$ and $r\leq q$ then we have $a^{(i,r)}=a^{(i,q)}$ for all integers $i<n.$ To show that the set $D_{n,a}\subset Q$ is dense, fix $q\in Q$ and by induction on $i\leq n+1$ build a decreasing sequence $q(0)\geq q(1)\geq \dots \geq q(i)\geq \dots$ so that 
\roster
\item $q=q(0)$
\item $q(i+1)\in D_{i,a}.$
\endroster
This is easily done, since at each step we have to meet only finitely many antichains. The above observation makes sure that by passing to stronger conditions we do not destroy the work done so far. The $q(n+1)\in D_{n,a}$ and $q(n+1)\leq q$ and the argument is complete. \qed C39
\enddemo

Now by the Proper Forcing Axiom there is a filter $H\subset Q$ meeting all the sets in the family $\{ D_{n,a}:n\in \omega,a\in[\oB]^{<\aA}\}.$ I define a function $g:H_\kk\to H_\kk$ by 
\roster
\item If $A\subset Q$ is a maximal antichain such that $H\cap A\neq 0$ then $g(A)=$ the unique element of $H\cap A.$
\item Otherwise the function $g$ is just identity.
\endroster

Let $M$ be the closure of $\oB$ under the functions $f,g.$ So $M\prec H_\kk,Q,\Delta\in M$ and $\oB\subset M.$ The following Claim completes the proof of the Lemma.

\proclaim
{Claim 40} Let $G=M\cap H.$ Then $G\subset Q\cap M$ is an $M$-generic filter.
\endproclaim

\demo
{Proof} For the genericity of $G,$ it is enough to prove that for any maximal antichain $A\subset Q$ in $M,$ We have $G\cap A\neq 0.$ So fix an antichain $A\in M.$ Since the model $M$ is chosen as a closure, there are a set $a\in[\oB]^{<\aA}$ and an integer $n$ such that $A$ belongs to the closure of $a$ under $f,g$ and is obtained after $n$ successive applications of the functions $f$ or $g.$ By the genericity of the filter $H\subset Q,$ there is a condition $q\in H\cap D_{n+1,a}.$ By the definition of the set $D_{n+1,a},$ the antichain $A$ belongs to the finite set $a^{(n,q)}$ and the condition $q$ has an element $x$ of $A$ above it. Since $q\leq x$ and $q\in H,$ we have $x\in H\cap A.$ Since the model $M$ is closed under the function $g,$ we have $x\in M$ and so $x\in G.$

We should verify that $G$ is a filter on $Q\cap M.$ Upwards closure follows from the same property of the filter $H.$ If $q$ and $r$ are two conditions in $G,$ then there is a lower bound of these two conditions in $H,$ but not {\it a priori} in $G.$ To remedy this defect we use the previous paragraph: by elementarity, in the model $M$ there is a maximal antichain $A\subset Q$ such that for $x\in A,$ either $r\perp x$ or $q\perp x$ or $x\leq q,r.$ By the above argument, there is $x\in A$ with $x\in G\subset H.$ But this $x$ must be compatible with both $q,r$ (since $H$ is a filter) and so it falls into the third category. Thus there is a lower bound of $q$ and $r$ in $G$ and $G$ is a filter. \qed C40,L38
\enddemo

\enddemo

The proof of Corollary 37 is now finished in the same fashion as the argument for Corollary 14.

\proclaim
{Main Theorem} If ZFC set theory is consistent then so is ZFC+``every complete Boolean algebra of uniform density $\aB$ contains a complete subalgebra isomorphic to $\CaB".$
\endproclaim

\demo
{Proof} The hard work has been done. The proof is now a routine iteration argument using Theorem 4 to deal with one algebra at a time. We give only an outline of the argument, since we believe that a reader that could bear with us up to here can easily provide the details. The scrupulous reader is advised to check with \cite {Scn} for every detail.

We start with a model of ZFC+GCH and set up a countable support iteration

$$\langle P_\aa:\aa\leq\oC,\dot Q_\aa:\aa<\oC\rangle$$

such that $P_\oC\Vdash$``every complete Boolean algebra of uniform density $\aB$ contains a complete subalgebra isomorphic to $\CaB".$ We shall have
\roster
\item the iterands are proper $\oC$-p.i.c. forcings of size $\aC.$
\endroster
Using a suitable bookkeeping device $\langle \tau _\aa:\aa\in\oC\rangle$ we shall browse through all potential $P_\oC$-names $\tau_\aa$ for separative posets of uniform density $\aB$ whose universe is $\oB.$ At all intermediate stages $\aa<\oC$ I shall have
\roster
\item"{(2)}" $P_\aa$ is a proper $\aC$-c.c. forcing notion of size $\aC$
\item"{(3)}" $P_\aa\Vdash$``GCH''
\endroster
These two properties hold true for any countable support iteration with property (1)--see Fact 32. So it will be possible, using Theorem 17 in $V^{P_\aa},$ to pick a $P_\aa$-name $\dot Q_\aa$ for a proper $\oC$-p.i.c. forcing of size $\aC$ so that
\roster
\item"{(4)}" $P_\aa\Vdash$``if $\tau_\aa$ is a separative poset of uniform density $\aB$ then $\dot Q\Vdash C_\aB\lessdot \tau_\aa".$
\endroster
For the final forcing $P_\oC,$ the following will be true:
\roster
\item"{(5)}" $P_\oC$ is a proper $\aC$-c.c. forcing--this holds by (1) and Fact 32.
\item"{(6)}" $|P_\oC|=\aC$--this is because the forcing $P_\oC$ is a direct limit of the forcings $P_\aa,\aa<\oC$ of size $\aC.$
\endroster
The properties (5),(6) make it possible to choose that suitable bookkeeping device $\langle \tau_\aa:\aa<\oC\rangle.$ Now by (5), $P_\oC$ does not collapse cardinals. We must verify that $P_\oC\Vdash$``every complete Boolean algebra of uniform density $\aB$ contains a complete subalgebra isomorphic to $\CaB".$ So let $\tau$ be a $P_\oC$-name such that $P_\oC\Vdash$``$\tau$ is a separative poset of uniform density $\aB$ with universe $\oB".$ Then for some $\aa<\oC$ we shall have that $\tau=\tau_\aa,$ $\tau$ is a $P_\aa$-name and $P_\aa\Vdash$``$\tau$ is a separative poset of uniform density $\aB$ with universe $\oB".$ By (4) we have $P_{\aa+1}\Vdash C_\aB\lessdot\tau.$ The statement $C_\aB\lessdot\tau$ is upwards absolute (Claim 15) and so we have even $P_\oC\Vdash C_\aB\lessdot\tau.$ As a result, $P_\oC\Vdash$`` the poset $C_\aB$ regularly embeds into every separative poset $\tau$ of uniform density $\aB"$ and the Theorem is proven.\qed MT
\enddemo

\subhead
{4. Towards higher densities}
\endsubhead

A natural question arises immediately upon seeing results \` a la Theorem 5: Is it possible to repeat such a feat for cardinalities higher than $\aB?$ We are very pessimistic about such a possibility; already the $\aC$ case seems to present unsurmountable difficulties. The following Theorem is the best negative result we can find in ZFC:

\proclaim
{Theorem 41} There is a separative partially ordered set $P$ of uniform density
 $\aleph _{\omega +1}$ such that $C_{\aleph _{\oA+1}}$ does not embed into it.
\endproclaim

It should be remarked that if e. g. the cardinals are the same in $V$ as in $L,$ the constructible universe, then I can find a poset as in Theorem 6 already in $L.$

\demo
{Proof} We shall need the following two facts from pcf theory.

\proclaim
{Lemma 42} Let $\langle \kappa _n:n\in \omega \rangle$ be an increasing
sequence of regular cardinals with $tcf(\prod _{n\in \omega }\kappa _n)$ mod fin $ =\lambda$ as witnessed by a  modulo finite increasing and cofinal sequence $\langle f_\beta :\beta<\lambda\rangle \subset
\prod _{n\in \omega }\kappa _n.$
Then there are ordinals $\beta _0<\beta _1<\lambda$ such that for all $n\in \omega$
we have $f_{\beta _0}(n)\leq f_{\beta _1}(n).$
\endproclaim

The proof is supplied below.

\proclaim
{Fact 43} (\cite {Sh2}) There is $\langle \kappa _n:
n\in \omega \rangle ,$ an increasing
sequence of regular cardinals $<\aleph _\omega$ with 
$tcf(\prod _{n\in \omega }\kappa _n)$ mod fin $=\aleph _{\omega +1}.$
\endproclaim

Now fix $\langle \kappa _n:n\in \omega \rangle ,$ an increasing
sequence of regular cardinals $<\aleph _\omega$ with 
$tcf(\prod _{n\in \omega }\kappa _n)$ mod fin $=\aleph _{\omega +1}$ and 
a modulo finite increasing and cofinal sequence $\langle f_\beta :\beta <\omega _{\omega +1} \rangle \subset\prod _{n\in \omega }\kappa _n.$ I am ready to define my partially ordered set $P:$

\proclaim
{Definition 44} The partially ordered set $P$ is a set of all pairs $\langle s,f
\rangle$ such that there are an integer $m$ with $s\in \prod _{n\in m }\kappa _n$ and an ordinal $\beta <\omega _{\omega +1}$ with $f=f_\beta \} .$

The order is defined by $\langle s^0,f^0\rangle \geq \langle s^1,f^1\rangle$
if $s^0\subset s^1,$ $\forall n\in dom(s^1)\setminus dom(s^0)$ $s^1(n)>f^0
(n)$ and $\forall n\notin dom(s^1)$ $f^1(n)\geq f^0(n).$ 
\endproclaim

\remark
{Explanation} So we add a function in $\prod _{n\in\oA}\kk_n$ which modulo finite dominates all the $f_\bb$'s. The $s$ part of a condition in $P$ is just a finite piece of this function.
\endremark

We prove now that the poset $C_{\aleph _{\oA+1}}$ does not embed into $P.$ Actually, more is true: if $c:\oA_{\oA+1}^V\to 2$ is a function in the generic extension by $P,$ then there is an infinite set $A\subset \oA_{\oA+1}$ in the ground model such that $c\restriction A$ is in the ground model again. Consequently, the function $c$ cannot be $C_{\aleph _{\oA+1}}$-generic over the ground model.

So let $p\in P,p\Vdash$``$\dot c:\omega _{\oA+1}^V\to 2$ is a function". We choose a sequence $\{ \langle s_\aa,f_{\bb_\aa}\rangle ,i_\aa:\aa\in\oA_{\oA+1}\}$ such that the following conditions are satisfied:

\roster
\item for each ordinal $\aa\in \oA_{\oA+1}$ I have $\langle s_\aa,f_{\bb_\aa}\rangle\in P,$ $i_\aa\in 2$
\item for each ordinal $\aa\in \oA_{\oA+1}$ I have $\langle s_\aa,f_{\bb_\aa}\rangle\Vdash$`` $\dot c(\aa)=i_\aa"$
\item for ordinals $\xi <\nu<\oA_{\oA+1}$ I have $\bb_\xi<\bb_\nu.$
\endroster

This is easily done. Now there are a set $S\subset \oA_{\oA+1}$ of full cardinality and a finite sequence $s$ such that for every ordinal $\aa\in S$ the constructed $s_\aa$ is just $s.$ We define the following partition $h$ of $S^2:$ for ordinals $\xi<\nu$ both in $S$ I set $h(\xi,\nu)=0$ if there is an integer $n$ such that $f_{\bb_\xi}(n)>f_{\bb_\nu}(n);$ otherwise, we let $h(\xi ,\nu)=1.$ By the Erd\H os-Dushnik-Miller theorem, we can have two cases:

\roster
\item There is a set $T\subset S$ of cardinality $\aleph _{\oA+1}$ homogeneous in $0.$ But this cannot happen since then the sequence $\langle f_{\bb_\aa}:\aa\in T\rangle\subset \prod _{n\in \oA}\kk_n$ contradicts Lemma 42. Notice that this sequence is indeed cofinal in $\prod _{n\in\oA}\kk_n$ since by (3) above, the set $\{ \bb_\aa:\aa\in T\}$ is cofinal in $\oA_{\oA+1}.$
\item There is a set $A\subset S$ of ordertype $\oA+1$ homogeneous in $1.$
\endroster

Since the first case leads to contradiction, the second case must happen. But then, if $\aa=max(A)$ and $\xi\in A,$ we have by the definition of the poset $P$ that $\langle s_\xi,f_{\bb_\xi}\rangle \geq \langle s_\aa,f_{\bb_\aa}\rangle.$ As a result, $\langle s_\aa,f_{\bb_\aa}\rangle\Vdash$``for every $\xi\in A,$ I have $\dot c(\xi)=i_\xi"$ and the argument is complete, since the condition $\langle s_\aa,f_{\bb_\aa}\rangle \leq p$ decides the values of $\dot c$ on an infinite set $A$ as desired. This leaves us with the last thing to demonstrate, namely Lemma 42.

\demo
{Proof of Lemma 42} The proof is quite technical and is modeled after Todorcevic's proof of a  similar fact about unbounded sequences of functions in $^\oA\oA$ \cite {T3}. Fix $\gamma <\lambda$ such that
$\{ s\in \bigcup _{m\in \omega }\prod _{n\in m}\kappa _n:\exists \beta
<\gamma$ $s\subset f_\beta \} =$
$\{ s\in \bigcup _{m\in \omega }\prod _{n\in m}\kappa _n:\exists \beta
<\lambda$ $s\subset f_\beta \} .$ This is possible since $\lambda >
sup\langle \kappa _n:n\in \omega \rangle$ is regular. I choose an integer $n_0$ and a set $S\subset \lambda$ of full cardinality so
that for every $n\geq n_0$ and for every $\beta \in S$ I have $f_\beta (n)\geq f_\gamma
(n).$ Define $T=
\{  s\in \bigcup _{m\in \omega }\prod _{n\in m}\kappa _n:$
$|\{ \beta \in S :s\subset f_\beta \} |=\lambda \} .$ So $T$ is a tree
of height $\omega .$ By induction on $n\in \omega$ simultaneously for all
$s\in T$ I define sets $A(s,n):$
\roster
\item $A(s,0)=\{ t\in T:s\subset t,lth(t)=lth(s)+1\}$
\item $A(s,n+1)=\{ t\in T:s\subset t,lth(t)=lth(s)+1,|A(t,n)|=
\kappa _{lth(t)}\} .$
\endroster

\proclaim
{Claim 45} There is $s\in T$ such that for all $n\in \omega$ $|A(s,n)|=
\kappa _{lth(s)}.$
\endproclaim

\demo
{Proof of the Claim} By contradiction. Assume that the Claim is false and
for any sequence $s\in T$ define $o(s)=min\{ n\in \omega :|A(s,n)|<\kappa _{lth(s)}\} .$
Choose $\delta \in \lambda$ such that for all $s\in \bigcup _{m\in \omega }
\prod _{n\in m}\kappa _n\setminus T$ we have
 $\{ \beta \in S:s\subset f_\beta \}
\subset \delta .$ We define a function $g \in \prod _{n\in \omega}\kappa _n$ by:

$$g(n)=max\{ f_\delta (n),sup \{ t(n):t\in \prod _{m\in n+1}\kappa _m
\cap T\text{ and } t\in A(t\restriction n,o(t\restriction n))\} \}$$

This is well-defined as the sets $A(s,o(s))$ are small. Now by the cofinality of the sequence $\langle
f_\beta :\beta \in S\rangle$ one can find an ordinal $\beta \in S$ and integer $n_1$
such that for all $n\geq n_1$ $f_\beta (n)\geq g(n).$ By our choice
of the ordinal $\delta$ we have that $f_\beta$ is a path through $T.$ It can be
easily verified now that the sequence of integers $\langle o(f _\beta\restriction n):n\geq n_1\rangle$
is strictly decreasing before it hits 0 for the first time. Let $n_2\geq
n_1$ be such that $o(f_\beta \restriction n_2)=0.$ So $|A(f_\beta \restriction
n_2,0)|<\kappa _{n_2}$ and since $f_\beta \restriction n_2+1\in 
A(f_\beta \restriction n_2,0)$ we obtain $f_\beta (n_2)<g(n_2),$
contradicting our choice of $n_1.$ \qed C45
\enddemo

To complete the proof of Lemma 42, choose a sequence $s\in T$ as in Claim 45. By my choice of $\gg,$ there is an ordinal
$\beta _0<\gamma$ such that $s\subset f_{\beta _0}.$ Since $f_{\bb_0}$ is modulo finite less than $f_\gg$ I can find an integer
$n_1\geq n_0$ such that $\forall n\geq n_1$ $f_{\beta _0}(n)\leq f_\gamma (n).$ Set $m=lth(s)$
and choose by induction finite sequences $s=s_m\subset s_{m+1}\subset \dots \subset s_{n_1}$
so that:
\roster
\item $s_j\in T,$ $lth(s_j)=j$
\item $s_{j+1}\in A(s_j,n_1-j),$ $s_{j+1}(j)\geq f_{\beta _0}(j).$
\endroster

This is possible since by induction on $j,$ $m\leq j \leq n_1$ one can verify
 that $|A(s_j,n_1-j)|=\kappa _j.$ Now pick $\beta _1\in S$ with 
$s_{n_1}\subset f_{\beta _1}.$ I claim that the ordinals $\bb_0<\bb_1$ exemplify the statement of the Lemma.

So I should show that for $n\in \omega,$ $f_{\beta _0}(n)\leq f_{\beta _1}(n).$ There are three cases. If $n<lth(s)$
then actually $f_{\beta _0}(n)= f_{\beta _1}(n).$ For $lth(s)\leq n<n_1$
the desired inequality follows from (2) above and for $n\geq n_1$ the inequality holds since
 $f_{\beta _0}(n)\leq f_\gamma (n)\leq f_{\beta _1}(n)$ (remember
$n\geq n_1\geq n_0).$ The argument is complete. \qed L42,T41
\enddemo
\enddemo

\subhead
{5. Open problems}
\endsubhead

There are several questions related to the Main Theorem left open in this thesis. The first two concern the structure of the real line in the resulting model.

\proclaim
{Problem 46} Assume that $\Bbb C(\aleph _1)$ embeds into every  algebra of uniform density $\aleph _1.$ Does it follow that $2^{\aleph _0}=\aleph _2?$
\endproclaim

\proclaim
{Problem 47} Assume that $\Bbb C(\aleph _1)$  embeds into every algebra of uniform density $\aleph _1.$ Does it follow that there is a Cohen real over $L?$
\endproclaim

Section 4 provides definite limitations for the possibility of obtaining results \` a la Theorem 5 for higher densities than $\aleph _1.$ In the positive direction we ask (motivated by \cite {FMS}):

\proclaim
{Problem 49} Is it consistent that $\Bbb C(\kappa )$ embeds into every separative partial order in $L$ of uniform density $\kappa?$ Is it implied by $0^\# ?$
\endproclaim

The following questions can hopefully inspire further development of my techniques for the $\aleph _1$ case:

\proclaim
{Problem 50} Is it consistent that every $\omega$-proper poset of size $\aleph _1$ is essentially c.c.c.? (I. e. there is a dense set $D\subset P$ of conditions such that for every $p\in D$ the poset $P\restriction p$ is c.c.c.)
\endproclaim

\proclaim
{Problem 51} Is it true that for any separative poset $P$ of uniform density $\kk$ which has a dense almost avoidable subset, we obtain $C_\kk\lessdot P?$
\endproclaim

\enddocument